\documentclass[a4paper, 10pt]{article}

\usepackage{amssymb}
\usepackage{amsmath}
\usepackage[dvips]{graphicx}% to include graphics
\usepackage[all]{xy}        % Xy-pic
\usepackage{alltt}          % print verbatim + LaTeX code
\usepackage{shortvrb}       % shorthand for \verb with \MakeShortVerb{} and \DeleteShortVerb{|}
\usepackage{amsthm}         % AMS theorem package
\usepackage{url}            % insert Internet addresses with \url
\usepackage{fancyhdr}       % fancy headers
\newcommand{\EG}{\ensuremath{\underline{E}G}}
\newcommand{\E}[1]{\ensuremath{\underline{E}#1}}   %enable the command to be called in both math & text mode
\newcommand{\F}{\ensuremath{\mathfrak{F}}}
\newcommand{\EF}[1]{\ensuremath{E_\mathfrak{F}#1}}
\newcommand{\orbitcat}{\ensuremath{\mathcal{O}_\mathfrak{F}G}}
\newcommand{\Fin}{\ensuremath{\mathfrak{Fin}}}
\newcommand{\Ind}{\textup{Ind}}

\newcommand{\Gmod}{\ensuremath{G\textup{-Mod}_{\mathfrak{F}}}}
\newcommand{\modG}{\ensuremath{\textup{Mod}_{\mathfrak{F}}\textrm{-}G}}
\newcommand{\stab}{\textup{stab}}
\newcommand{\ttimes}{\!\times\!}

\newcommand{\modGH}{\textup{Mod}_{\mathfrak{F} \times \mathfrak{F'}}\textup{-}(G \ttimes H)}

\newcommand{\Z}{\ensuremath{\mathbb{Z}}}
\newcommand{\R}{\ensuremath{\mathbb{R}}}
\newcommand{\RR}{\ensuremath{\mathcal{R}}}

\newtheorem{proposition}{Proposition}[section]
\newtheorem{theorem}[proposition]{Theorem}
\newtheorem{corollary}[proposition]{Corollary}
\newtheorem{lemma}[proposition]{Lemma}
\newtheorem{defn}[proposition]{Definition}
\theoremstyle{remark} \newtheorem{remark}{Remark}%

\numberwithin{equation}{section} %Change equation numbering into section.number

\title{Equivariant $K$-homology for some Coxeter groups}
\author{Rub\'en S\'anchez-Garc\'\i{}a\thanks{Funded by the EPSRC and the School of
Mathematics, University of Southampton}}
\date{\today}

\begin{document}

\maketitle

\begin{abstract}
    \noindent We obtain the equivariant $K$-homology of the
    classifying space $\underline{E}W$ for $W$ a right-angled or,
    more generally, an even Coxeter group. The
    key result is a formula for the relative Bredon homology of
    $\underline{E}W$ in terms of Coxeter cells. Our
    calculations amount to the $K$-theory of the reduced
    $C^*$-algebra of $W$, via the Baum-Connes assembly map.
\end{abstract}

\section{Introduction}
Consider a discrete group $G$. The Baum-Connes conjecture
\cite{BCH94} identifies the $K$-theory of the reduced
$C^*$-algebra of $G$, $C^*_r(G)$, with the equivariant
$K$-homology of a certain classifying space associated to $G$.
This space is called the classifying space for proper actions,
written \EG{}. The conjecture states that a particular map between
these two objects, called the assembly map,
    $$
    \mu_i : K_i^G(\E G) \longrightarrow K_i(C_r^*(G))
    \quad i \ge 0
    \, ,
    $$
is an isomorphism. Here the left hand side is the equivariant
$K$-homology of \EG{} and the right hand side is the $K$-theory of
$C^*_r(G)$. The conjecture can be stated more generally
\cite[Conjecture 3.15]{BCH94}.

The equivariant $K$-homology and the assembly map are usually
defined in terms of Kasparov's $KK$-theory. For a discrete group
$G$, however, there is a more topological description due to Davis
and L\"uck \cite{DL98}, and Joachim \cite{Joa03} in terms of
spectra over the orbit category of $G$. We will keep in mind this
topological viewpoint (cf.~Mislin's notes in
\cite{MislinValette:ProperGroupActions}).

Part of the importance of this conjecture is due to the fact that
it is related to many other relevant conjectures \cite[\S
7]{MislinValette:ProperGroupActions}. Nevertheless, the conjecture
itself allows the computation of the $K$-theory of $C^*_r(G)$ from
the $K^G$-homology of \EG{}. In turn, this $K$-homology can be
achieved by means of the Bredon homology of \EG, as we explain
later.

In this article we focus our attention on (finitely generated)
Coxeter groups. These groups are well-known in geometric group
theory as groups generated by reflections (elements of order 2)
and only subject to relations in the form $(st)^n = 1$. Coxeter
groups have the Haagerup property \cite{BJS88} and therefore
satisfy the Baum-Connes conjecture \cite{HK97}.

For $W$ a Coxeter group, we consider a model of the classifying
space $\underline{E}W$ called the Davis complex. We obtain a
formula for the relative Bredon homology of this space in terms of
Coxeter cells (Theorem \ref{thm:FormulaCoxeterCells}). From this,
we can deduce the Bredon homology of $\underline{E}W$ in same
cases; for instance, for right-angled Coxeter groups and, more
generally, even Coxeter groups (Theorems \ref{thm:BredonHomRA} and
\ref{thm:BredonHomEven}). The equivariant $K$-homology of
$\underline{E}W$ follows immediately. Since Coxeter groups satisfy
Baum-Connes, our results also amount to the $K$-theory of the
corresponding reduced $C^*$-algebra.

These results appear in the author's PhD thesis \cite[Chapter
5]{SanchezGarciaPhD}. I would like to thank my PhD supervisor Ian
Leary for his guidance through this research project, and
particularly for suggesting the use of the Coxeter cell structure
and relative Bredon homology.

\section{Preliminaries}\label{section:Preliminaries}
\subsection{Classifying spaces}\label{section:ClassifSpaces}
The classifying space \EG{} that appears in the Baum-Connes
conjecture is a particular case of a more general construction.

Let $G$ be a discrete group. A $G$-CW-complex is a CW-complex with
a $G$-action permuting the cells and such that if a cell is sent
to itself, it is done by the identity map. Let \F{} be a non-empty
family of subgroups of $G$ closed under conjugation and passing to
subgroups. A \emph{model for \EF{G}} is a $G$-CW-complex $X$ such
that (1) all cell stabilizers are in \F{}; (2) for any
$G$-CW-complex $Y$ with all cell stabilizers in \F, there is a
$G$-map $Y \rightarrow X$, unique up to $G$-homotopy equivalence.
The last condition is equivalent to the statement that for each $H
\in \F$, the fixed point subcomplex $X^H$ is contractible.

For the family with just the trivial subgroup, we obtain $EG$, a
contractible free $G$-CW-complex whose quotient $BG$ is the
classifying space for principal $G$-bundles; when $\F =
\mathfrak{Fin}(G)$, the family of all finite subgroups of $G$,
this is the definition of \EG{}.

It can be shown that general classifying spaces $\EF{G}$ always
exists. They are clearly unique up to $G$-homotopy. See \cite[\S
2]{BCH94} or \cite{MislinValette:ProperGroupActions} for more
information and examples of \EG.

\subsection{Bredon (co)homology}\label{section:BredonHomology}
Given a group $G$ and a family \F{} of subgroups, we will write
\orbitcat{} for the \emph{orbit category}. The objects are left
cosets $G/K$, $K \in \F$, and morphisms the $G$-maps $\phi: G/K
\rightarrow G/L$. Such a $G$-map is uniquely determined by its
image $\phi(K) = gL$, and we have $g^{-1}Kg \subset L$.
Conversely, such $g \in G$ defines a $G$-map, which will be
written $R_g$.

A covariant (resp.~contravariant) functor $M \colon \orbitcat
\rightarrow \mathfrak{Ab}$ is called a \emph{left}
(resp.~\emph{right}) \emph{Bredon module}. The category of left
(resp.~right) Bredon modules and natural transformations is
written \Gmod{} (resp.~\modG). It is an abelian category, and we
can use homological algebra to define Bredon homology (see
\cite[pp.~7-10]{MislinValette:ProperGroupActions}). Nevertheless,
we give now a practical definition.

Consider a $G$-CW-complex $X$, and write $\mathfrak{Iso}(X)$ for
the family of isotropy subgroups $\{ \stab(x), x \in X\}$. Let
$\F$ be a family of subgroups of $G$ containing
$\mathfrak{Iso}(X)$, and $M$ a left Bredon module. The
\emph{Bredon homology groups} $H^\F_i \left(X;M\right)$ can be
obtained as the homology of the following chain complex $(C_*,
\partial_*)$.
Let $\{e_\alpha\}$ be orbit representatives of the $d$-cells ($d
\ge 0$) and write $S_\alpha$ for $\stab(e_\alpha) \in \F$. Define
$$
    C_d = \bigoplus_\alpha M \left( G/S_\alpha \right) \, .
$$
If $g \cdot e'$ is a typical $(d-1)$-cell in the boundary of
$e_\alpha$ then $g^{-1} \cdot \stab(e_\alpha) \cdot g \subset
\stab(e')$, giving a $G$-map (write $S'$ for $\stab(e')$)
$$
    R_g: G/S_\alpha \rightarrow G/S'\, ,
$$
which induces a homomorphism $M(\phi)\colon M \left(
G/S_\alpha\right) \rightarrow M\left(G/S'\right)$, usually written
$(R_g)_*$. This yields a differential $\partial_d\colon C_d
\rightarrow C_{d-1}$, and the Bredon homology groups
$H^\F_i\left(X;M\right)$ correspond to the homology of $(C_*,
\partial_*)$. Observe that the definition is independent of the
family \F{} as long as it contains $\mathfrak{Iso}(X)$.

Bredon cohomology is defined analogously, for $M$ a right Bredon
mo\-du\-le (a contravariant $M$ will reverse the arrows $(R_g)^* =
M(\phi)\colon M(G/S') \rightarrow M(G/S_\alpha)$ so that
$\partial: C_{d-1} \rightarrow C_d$).

The Bredon homology of a group $G$ with coefficients in $M \in
\Gmod$ can be defined in terms of a Tor functor
(\cite[Def.~3.12]{MislinValette:ProperGroupActions}). If \F{} is
closed under conjugations and taking subgroups then
$$ H_i^\F(G;M) \cong H_i^\F(E_\F{}G;M)\, ,$$
which may as well be taken as a definition.

We are interested in the case $X = \EG$, $\F = \Fin(G)$ and $M =
\mathcal{R}$ the complex representation ring, considered as a
Bredon module as follows. On objects we set
    $$
        \mathcal{R}(G/K) = R_\mathbb{C}(K), \quad K \in \Fin(G)
    $$
the ring of complex representations of the finite group $K$
(viewed just as an abelian group), and for a $G$-map $R_g: G/K
\rightarrow G/L$, we have $g^{-1}Kg \subset L$ so that $(R_g)_*:
R_\mathbb{C}(K) \rightarrow R_\mathbb{C}(L)$ is given by induction
from $g^{-1}Kg$ to $L$ after identifying $R_\mathbb{C}(g^{-1}Kg)
\cong R_\mathbb{C}(K)$.

\subsection{Equivariant $K$-homology}\label{section:EquivariantKhomology}
There is an equivariant version of $K$-homology, denoted
$K^G_i(-)$ and defined in \cite{DL98} (see also \cite{Joa03})
using spaces and spectra over the orbit category of $G$. It was
originally defined in \cite{BCH94} using Kasparov's $KK$-theory.
We will only recall the properties we need.

Equivariant $K$-theory satisfies Bott mod-2 periodicity, so we
only consider $K^G_0$ and $K^G_1$. For any subgroup $H$ of $G$, we
have
$$
    K^G_i\left(G/H\right) = K_i\left(C^*_r(H)\right) \, ,
$$
that is, its value at one-orbit spaces corresponds to the
$K$-theory of the reduced $C^*$-algebra of the typical stabilizer.
If $H$ is a \emph{finite} subgroup then $C^*_r(H) = \mathbb{C}H$
and
$$
    K_i^G\left(G/H\right) = K_i(\mathbb{C}H) = \left\{
    \begin{array}{cl} R_\mathbb{C}(H) & i = 0\, ,\\ 0 & i = 1 \, .
    \end{array} \right.
$$
This allows us to view $K_i^G(-)$ as a Bredon module over
$\mathcal{O}_\Fin G$.

We can use an equivariant Atiyah-Hirzebruch spectral sequence to
compute the $K^G$-homology of a proper $G$-CW-complex $X$ from its
Bredon homology (see
\cite[pp.~49-50]{MislinValette:ProperGroupActions} for details),
as
$$
    E^2_{p,q} = H^\Fin_p\left( X; K^G_q(-) \right) \Rightarrow
    K^G_{p+q}\left(X\right)\, .
$$
In the simple case when Bredon homology concentrates at low
degree, it coincides with the equivariant $K$-homology:
\begin{proposition}\label{prop:LowDegreeBredonHomology}
    Write $H_i$ for $H_i^\Fin(X;\RR)$ and $K_i^G$ for $K_i^G(X)$.
    If $H_i = 0$ for $i \ge 2$ then $K^G_0 = H_0$ and $K^G_1 = H_1$.
\end{proposition}
\begin{proof} The Atiyah-Hirzebruch spectral sequence collapses at
the 2-page.
\end{proof}

%%% NEW SECTION %%%%%%%%%%%%%%%%%%%%%%%%%%%%%%%%%%%%%%%%%%%%%%%%%%%%
\section{K\"unneth formulas for Bredon homology}\label{section:KunnethFormulas}
We will need K\"unneth formulas for (relative) Bredon homology. We
devote this section to state such formulas, for the product of
spaces and the direct product of groups. Theorem
\ref{thm:KunnethFormulaXY} holds in more generality (see
\cite{GG99}) but we state the result we need and sketch a direct
proof; details can be found in \cite[Chapter 3]{SanchezGarciaPhD}.
Results similar to those explained here are treated in
\cite{LeonardiPhD_inprep}.

\subsection*{K\"unneth formula for $X \ttimes Y$}
Let $X$ be a $G$-CW-complex and $Y$ a $H$-CW-complex. Let \F{}
(resp.~$\F'$) be a family of subgroups of $G$ (resp.~of $H$)
containing $\mathfrak{Iso}(X)$ (resp.~$\mathfrak{Iso}(Y)$). Then
$X\ttimes Y$ is a $(G\ttimes H)$-CW-complex (with the compactly
generated topology) and
$$ \mathfrak{Iso}(X\ttimes Y) = \mathfrak{Iso}(X) \times
\mathfrak{Iso}(Y) \subset \F \times \F'\, .$$%
Moreover, the orbit category $\mathcal{O}_{\mathfrak{F} \times
\mathfrak{F'}}(G \ttimes H)$ is isomorphic to
$\mathcal{O}_\mathfrak{F}G \times \mathcal{O}_\mathfrak{F'}H$.

Given $M \in \Gmod$, $N \in H\textup{-Mod}_{\mathfrak{F'}}$ we
define their tensor product over $\mathbb{Z}$ as the composition
of the two functors
$$ \xymatrix{
M \otimes N: \mathcal{O}_\mathfrak{F}G \times
\mathcal{O}_\mathfrak{F'}H \ar[r]^-{M\times N} & \mathfrak{Ab}
\ttimes \mathfrak{Ab} \ar[r]^-{\otimes_\mathbb{Z}} &
\mathfrak{Ab}}
$$
considered a Bredon module over $\mathcal{O}_{\mathfrak{F} \times
\mathfrak{F'}}(G \ttimes H)$ $\cong$ $\mathcal{O}_\mathfrak{F}G
\times \mathcal{O}_\mathfrak{F'}H$. We can easily extend this
tensor product to chain complexes of Bredon modules.

\begin{theorem}[K\"unneth Formula for Bredon Homology]\label{thm:KunnethFormulaXY}
Consider $X$, $Y$, $\F$, $\F'$, $M$, $N$ as above, with the
property that $M(G/K)$ and $N(H/K')$ are free for all $K \in
\mathfrak{Iso}(X)$, $K' \in \mathfrak{Iso}(Y)$. Then for every $n
\geq 0$, there is a split exact sequence
\begin{eqnarray*}
    0 \rightarrow \bigoplus_{i+j=n} \Big( H^{\mathfrak{F}}_i (X;M) \otimes_\mathbb{Z}
    H^{\mathfrak{F'}}_j (Y;N) \Big) \rightarrow H_n^{\mathfrak{F} \times \mathfrak{F'}}
    (X \ttimes Y; M \otimes N) \\
     \rightarrow \bigoplus_{i+j=n-1} \textup{Tor} \Big( H^{\mathfrak{F}}_i
    (X;M),H^{\mathfrak{F'}}_j (Y;N)\Big) \rightarrow 0
\end{eqnarray*}
\end{theorem}
\begin{proof}\textit{(sketch)} Consider the chain complexes of abelian groups
\begin{eqnarray*}
    D_* = \underline{C_*(X)} \otimes_{\mathfrak{F}} M \quad \textrm{and} \quad
    D'_* = \underline{C_*(Y)} \otimes_{\mathfrak{F'}} N\, ,
\end{eqnarray*}
where $\underline{C_*(-)}$ is the chain complex of Bredon modules
defined in \cite[p.~10-11]{MislinValette:ProperGroupActions}, and
$-\otimes_\F-$ is the tensor product defined in
\cite[p.~14]{MislinValette:ProperGroupActions}. Then, by
Definition 3.13 in \cite{MislinValette:ProperGroupActions},
$$
    H_i (D_*) = H^{\mathfrak{F}}_i (X;M)\, , \quad H_j (D'_*)=
    H^{\mathfrak{F'}}_j (Y;N)\, .
$$
We now use the following lemmas, which follow from the definitions
of $\underline{C_*(-)}$ and $-\otimes_\F-$.

\begin{lemma}\label{proposition:12}
Given $\underline{C_*}$ a chain complex in $\modG$,
$\underline{C'_*}$ a chain complex in
$\textrm{Mod}_{\mathfrak{F'}}\textrm{-}H$ and Bredon modules $M
\in \Gmod$, $N \in H\textrm{-Mod}_{\mathfrak{F'}}$, the following
chain complexes of abelian groups are isomorphic
\begin{displaymath}
    \Big( \underline{C_*} \otimes_{\mathfrak{F}} M \Big) \otimes \Big( \underline{C'_*}
    \otimes_{\mathfrak{F'}} N \Big) \cong
    \Big( \underline{C_*} \otimes \underline{C'_*} \Big)
    \otimes_{\mathfrak{F}\times\mathfrak{F'}} \Big( M \otimes N \Big) .
\end{displaymath}
\end{lemma}
\begin{lemma}\label{proposition:11} The following chain
complexes in $\modGH$ are isomorphic
\begin{displaymath}
    \underline{C_*(X \ttimes Y)} \cong \underline{C_*(X)} \otimes
    \underline{C_*(Y)}\ .
\end{displaymath}
\end{lemma}
\noindent By Lemma \ref{proposition:12},
$$
D_* \otimes_\Z D_*' \cong \left(\underline{C_*(X )} \otimes
\underline{C_*(Y)}\right) \otimes_{\mathfrak{F}\times
\mathfrak{F'}} \left(M \otimes N\right)
$$
and therefore, by Lemma \ref{proposition:11},
$$
    H_n(D_* \otimes_\mathbb{Z} D_*') = H_n^{\mathfrak{F} \times \mathfrak{F'}}
    (X \ttimes Y; M \otimes N)\, .
$$
The result is now a consequence of the ordinary K\"unneth formula
of chain complexes of abelian groups. Observe that the chain
complexes $D_*$ and $D'_*$ are free since
\begin{displaymath}
    D_n = \underline{C_n(X)} \otimes_{\mathfrak{F}} M
    \cong \bigoplus_\alpha M(G/S_\alpha)\,
\end{displaymath}
(see \cite[p.~14-15]{MislinValette:ProperGroupActions}), and
similarly for $D'_n$.
\end{proof}
There is also a K\"unneth formula for \emph{relative} Bredon
homology (see \cite[p.~12]{MislinValette:ProperGroupActions} for
definitions).
\begin{theorem}[K\"unneth Formula for Relative Bredon Homology]\label{thm:KunnethFormulaRelative}
Under the same hypothesis of Theorem \ref{thm:KunnethFormulaXY},
plus $A \subset X$ a $G$-CW-subcomplex and $B \subset Y$ a
$H$-CW-subcomplex, we have that for every $n \geq 0$ there is a
split exact sequence
\begin{eqnarray*}
    \lefteqn{0 \longrightarrow \bigoplus_{i+j=n} \Big( H^{\mathfrak{F}}_i (X,A;M) \otimes_\mathbb{Z}
    H^{\mathfrak{F'}}_j (Y,B;N) \Big) } \\[0.5em]
    & & \longrightarrow H_n^{\mathfrak{F} \times \mathfrak{F'}}
    (X \ttimes Y, A \ttimes Y \cup X \ttimes B; M \otimes N)
    \\[0.6em]
     & & \longrightarrow\bigoplus_{i+j=n-1} \textup{Tor} \Big( H^{\mathfrak{F}}_i
    (X,A;M),H^{\mathfrak{F'}}_j (Y,B;N)\Big) \longrightarrow
    0 \ .
\end{eqnarray*}
\end{theorem}
\noindent The proof is the same as the non-relative formula but
using the chain complexes
$$
    \begin{array}{c}
        D_* = \underline{C_*(X,A)} \otimes_{\mathfrak{F}} M\\
        D'_* = \underline{C_*(Y,B)} \otimes_{\mathfrak{F'}} M
    \end{array}
$$
and a relative version of Lemma \ref{proposition:11}.

\subsection*{K\"unneth formula for $G \ttimes H$}
We will work in general with respect to a class of groups \F{}
(formally, \F{} is a collection of groups closed under
isomorphism, or a collection of isomorphism types of groups). Our
main example will be $\mathfrak{Fin}$, the class of finite groups.

Given a specific group $G$, let us denote by $\F(G)$ the family of
subgroups of $G$ which are in \F{}, and by $H^\F_*(-)$ the Bredon
homology with respect to $\F(G)$. Suppose now that \F{} is closed
under taking subgroups (note that it is always closed under
conjugation). Then so is $\F(G)$ and we will write $\EF{}G$ for
$E_{\F(G)}G$.

We want to know when  $E_\F{}G \times E_\F{}H$ is a model for
$E_\F(G\times H)$.
\begin{proposition} \label{prop:E_F(GxH)}
Let \F{} be a class of groups which is closed under subgroups,
finite direct products and homomorphic images. Then, for any
groups $G$ and $H$, and $X$ and $Y$ models for $\EF{}G$
respectively $\EF{}H$, the space $X \times Y$ is a model for
$\EF(G\ttimes H)$.
\end{proposition}
\begin{proof} If $(x,y) \in X \times Y$ then
$$\stab_{G\times H}(x,y) = \stab_G(x) \times \stab_H(y) \in
\F(G) \times \F(H) \subset \F(G\times H)\, .$$ Consider the
projections $\pi_1: G \ttimes H \rightarrow G$ and $\pi_2: G
\ttimes H \rightarrow H$. If $K \in \F(G \ttimes H)$,
$$ (X \times Y)^K = X^{\pi_1(K)} \times Y^{\pi_2(K)} $$
contractible since $\pi_1(K) \in \F(G)$ and $\pi_2(K) \in \F(H)$.
\end{proof}

\begin{remark}
The proposition is not true if we remove any of the two extra
conditions on \F.
\end{remark}

\begin{remark} Not every family of subgroups of $G$ can be
written as $\F(G)$ (for instance, $\F(G)$ is always closed under
conjugation). To be more general, one can prove that given two
families $\F_1$ and $\F_2$ of subgroups of $G$ respectively $H$,
$E_{\F_1}G \times E_{\F_2}H$ is a model for $\EF(G\times H)$,
where \F{} is the smallest family closed under subgroups and
containing $\F_1 \times \F_2$. Nevertheless, we are interested in
group theoretic properties (as `being finite'), so we think of
$\F$ as the class of groups with the required property.
\end{remark}

Thus we can apply the K\"unneth formula (Theorem
\ref{thm:KunnethFormulaXY}) to $E_\F(G\ttimes H)$. We obtain
Bredon homology groups with respect to the family $\F(G) \times
\F(H)$ instead of $\F(G\ttimes H)$, but both families contain the
isotropy groups of $E_\F{}G \times E_\F{}H$, so the Bredon
homology groups are the same (cf.~Section
\ref{section:BredonHomology}).
\begin{theorem}[K\"unneth formula for $G\ttimes H$]\label{thm:KunnethFormulaGH}
Let \F{} be a class of groups closed under taking subgroups,
direct product and homomorphic images. For every $n \geq 0$ there
is a split exact sequence
        \begin{eqnarray*}
        0 \rightarrow \bigoplus_{i+j=n} \Big( H^{\mathfrak{F}}_i (G;M) \otimes_\mathbb{Z}
        H^{\mathfrak{F}}_j (H;N) \Big) \rightarrow
        H_n^{\mathfrak{F}}
        (G \ttimes H; M \otimes N) \rightarrow \\
        \bigoplus_{i+j=n-1} \textup{Tor} \Big(
        H^{\mathfrak{F}}_i
        (G;M),H^{\mathfrak{F}}_j (H;N)\Big) \rightarrow 0\, .
        \end{eqnarray*}
\end{theorem}
\begin{remark}
The analogous result for relative Bredon homology also holds, by
Theorem \ref{thm:KunnethFormulaRelative}.
\end{remark}

\subsection*{Application to proper actions and coefficients in the representation ring}
The original motivation for stating a K\"unneth formula was the
Bredon homology of $G \ttimes H$ for proper actions (i.e.~$\F =
\mathfrak{Fin}$) and coefficients in the representation ring. We
only need a result on the coefficients.

For finite groups $P$ and $Q$, the representation ring
$R_\mathbb{C}(P \ttimes Q)$ is isomorphic to $R_\mathbb{C}(P)
\otimes R_\mathbb{C}(Q)$. The same is true at the level of Bredon
coefficient systems. Given groups $G$ and $H$, consider the left
Bredon modules
$$ \mathcal{R}^G \in G\textrm{-Mod}_{\mathfrak{Fin}(G)}\; ,\ \mathcal{R}^H \in
    H\textrm{-Mod}_{\mathfrak{Fin}(H)} \;\textrm{ and }\; \mathcal{R}^{G \ttimes H}
    \in (G \ttimes H)\textrm{-Mod}_{\mathfrak{Fin}(G \ttimes
    H)}\, .
$$
Let us denote by $\widetilde{\mathcal{R}}^{G \ttimes H}$ the
restriction of $\mathcal{R}^{G \ttimes H}$ to $\mathfrak{Fin}(G)
\ttimes \mathfrak{Fin}(H)$.

\begin{proposition}\label{proposition:RepresentationRingGxH}
These two Bredon modules are naturally isomorphic
\begin{displaymath}
    \mathcal{R}^G \otimes \mathcal{R}^H \cong \widetilde{\mathcal{R}}^{G \ttimes H}\, .
\end{displaymath}
\end{proposition}
\begin{proof} Consider the isomorphism $\Theta_{P,Q} : R_\mathbb{C}(P \ttimes Q) \cong
R_\mathbb{C}(P) \otimes R_\mathbb{C}(Q)$ and the universal
property of the tensor product to obtain a map $\alpha$ as
\begin{displaymath}
    \xymatrix{
    R_\mathbb{C}(P) \ttimes R_\mathbb{C}(Q)  \ar[r] \ar[dr]_\alpha &
    R_\mathbb{C}(P) \otimes R_\mathbb{C}(Q) \ar[d]_\cong^{\Theta_{P,Q}}\\
      & R_\mathbb{C}(P \ttimes Q)\, .}
\end{displaymath}
That is, $\Theta_{P,Q}(\rho \otimes \tau) = \alpha(\rho, \tau)$,
where the latter denotes the ordinary tensor product of
representations and $\rho \otimes \tau$ is their formal tensor
product in $R_\mathbb{C}(P) \otimes R_\mathbb{C}(Q)$. Extend this
map to a natural isomorphism $\Theta: \mathcal{R}^G \otimes
\mathcal{R}^H \rightarrow \widetilde{\mathcal{R}}^{G \ttimes H}$.
It only remains to check naturality, that is, the commutativity of
\begin{equation*}
    \xymatrix{
    R_\mathbb{C}(P) \otimes R_\mathbb{C}(Q)
    \ar[r]^{\Theta_{P,Q}} \ar[d]_{(R_g)_* \otimes (R_h)_*} & R_\mathbb{C}(P \times Q)
    \ar[d]^{(R_{(g,h)})_*} \\
    R_\mathbb{C}(P') \otimes R_\mathbb{C}(Q')
    \ar[r]^{\Theta_{P',Q'}} & R_\mathbb{C}(P' \times Q')\, ,}
\end{equation*}
for any $(G\times H)$-map $R_{(g,h)} = R_g \times R_h$. To do
this, one can compute both sides of the square and show that the
two representations have the same character, using the following
two observations. Denote by $\chi(-)$ the character of a
representation. Since $\alpha(-,-)$ is the tensor product of
representations we have $\chi\left(\alpha(\tau^1,\tau^2)\right) =
\chi(\tau^1) \cdot \chi(\tau^2)$. If $R_g : G/K \rightarrow G/L$
is a $G$-map, $\rho \in R_\mathbb{C}(K)$ and $t \in G$ then
\begin{equation*}\label{eqn:LemmaCharIndRepr}
    \chi \left( \left( R_{g}\right)_* \rho \right)
    (t) = \frac{1}{|K|} \sum_* \chi(\rho) (k) \, ,
\end{equation*}
where the sum $*$ is taken over the $s \in G$ such that $k =
gsts^{-1}g^{-1} \in K$. We leave the details to the reader.
\end{proof}

For an arbitrary group $G$, write $H_i^{\mathfrak{Fin}}(G;
\mathcal{R})$ for $H_i^{\mathfrak{Fin}}(G; \mathcal{R}^G)$.
    \begin{corollary} \label{cor:KunnethFormulaGxH}
    For every $n \geq 0$ there is a split exact sequence
        \begin{eqnarray*}
        0 \rightarrow \bigoplus_{i+j=n} \Big( H^{\mathfrak{Fin}}_i (G;\mathcal{R}) \otimes_\mathbb{Z}
        H^{\mathfrak{Fin}}_j (H;\mathcal{R}) \Big) \rightarrow
        H_n^{\mathfrak{Fin}}
        (G \ttimes H; \mathcal{R}) \rightarrow \\
        \bigoplus_{i+j=n-1} \textup{Tor} \Big(
        H^{\mathfrak{Fin}}_i
        (G;\mathcal{R}),H^{\mathfrak{Fin}}_j (H;\mathcal{R})\Big) \rightarrow
        0 \, .
        \end{eqnarray*}
    \end{corollary}
\begin{remark}
Again, there is a similar statement for relative Bredon homology.
\end{remark}

\section{Coxeter groups and the Davis complex}
We briefly recall the definitions and basic properties of Coxeter
groups, and describe a model of the classifying space for proper
actions. Most of the material is well-known and can be found in
any book on the subject, as \cite{Humphreys:CoxeterGroups}.

\subsection{Coxeter groups}
Suppose $W$ is a group, $S = \{ s_1, \ldots , s_N \}$ a finite
subset of elements of order 2 which generate $W$. Let $m_{ij}$
denote the order of $s_is_j$. Then $2 \le m_{ij} = m_{ji} \le
\infty$ if $i \neq j$ and $m_{ii} = 1$. We call the pair $(W,S)$ a
\emph{Coxeter system} and $W$ a \emph{Coxeter group} if $W$ admits
a presentation
$$\label{def:CoxeterGroup}
        \langle\, s_1, \ldots , s_N \;|\; (s_is_j)^{m_{ij}} = 1\, \rangle \ .
$$
We call $N$ the \emph{rank} of $(W,S)$. The \emph{Coxeter diagram}
is the graph with vertex set $S$ and one edge joining each pair
$\{s_i, s_j\}$, $i \neq j$, with label $m_{ij}$, and the
conventions: if $m_{ij} = 2$ we omit the edge; if $m_{ij} = 3$ we
omit the label. A Coxeter system is \emph{irreducible} if its
Coxeter diagram is connected. Every Coxeter group can be
decomposed as a direct product of irreducible ones, corresponding
to the connected components of its Coxeter diagram. All
\emph{finite} irreducible Coxeter systems have been classified; a
list of their Coxeter diagrams can be found, for instance, in
\cite[p.~32]{Humphreys:CoxeterGroups}.

Every Coxeter group can be realized as a group generated by
reflections in $\R^N$. Namely, there is a faithful \emph{canonical
representation} $W \rightarrow GL_N(\R)$ (see \cite[\S
5.3]{Humphreys:CoxeterGroups} for details). The generators in $S$
correspond to reflections with respect to hyperplanes. These
hyperplanes bound a \emph{chamber}, which is a strict fundamental
domain (examples below).

\subsection*{Coxeter cells}
Suppose now that $W$ is finite. Take a point $x$ in the interior
of the chamber. The orbit $Wx$ is a finite set of points in
$\R^N$. Define the \emph{Coxeter cell} $C_W$ as the convex hull of
$Wx$. It is a complex polytope by definition.\\[10pt]
\noindent\textbf{Examples:}
\begin{enumerate}
\item[(1)] Rank 0: $W = \{ 1 \}$ the trivial group and $C_W$ is a
point.

\item[(2)] Rank 1: $W = C_2$ the cyclic group of order 2,
$\mathbb{R}^N = \mathbb{R}$, the $W$-action consists on reflecting
about the origin, and there are two chambers, $(-\infty, 0]$ and
$[0,+\infty)$. Take a point $x$ in the interior of, say,
$[0,+\infty)$; then $Wx = \{ -x, x\}$ and $C_W = [-x,x]$, an
interval.

\item[(3)] Rank 2: $W = D_{n}$ dihedral group of order $2n$,
$\mathbb{R}^N = \mathbb{R}^2$ and the canonical representation
identifies $W$ with the group generated by two reflections respect
to lines through the origin and mutual angle $\pi / n$. A chamber
is the section between the two lines and the $W$-orbit $Wx$ is the
vertex set of a $2n$-gon (see Figure \ref{fig:hexagon}).
Therefore, $C_W$ is a $2n$-gon.

\begin{figure}[hbt]
   \begin{center}
   \includegraphics[scale = 0.3]{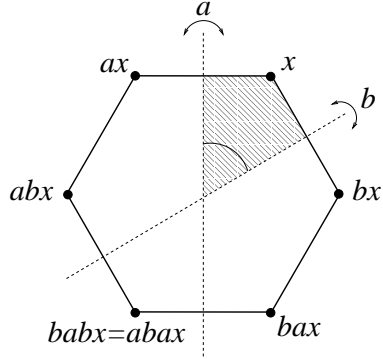}
   \caption{Coxeter cell for a dihedral group.}
   \label{fig:hexagon}
   \end{center}
   \end{figure}

\item[(4)] \label{eg:CoxeterCellProduct}If $W = W_1 \times W_2$
then $C_W = C_{W_1} \times C_{W_2}$.
\end{enumerate}

\subsection*{The Davis complex}
It is a standard model of $\underline{E}W$ introduced by M. Davis
in \cite{Davis83}. Given $T \subset S$, the group $W_T = \langle T
\rangle$ generated by $T$ is called a \emph{special subgroup}. It
can be shown that $(W_T,T)$ is a Coxeter system, and that $W_T
\cap W_{T'} = W_{T \cap T'}$. The subset $T$ is called
\emph{spherical} is $W_T$ is finite. We will write $\mathcal{S}$
for the set $\{ T \subset S \, | \, T \textrm{ is spherical} \}$.

A \emph{spherical coset} is a coset of a finite special subgroup,
that is, $wW_T$ for some $w \in W$, $T \in \mathcal{S}$. Note that
$wW_T = w'W_{T'}$ if and only if $T = T'$ and $w^{-1}w' \in W_T$
(use the so-called \emph{deletion condition}
\cite[\S{}5.8]{Humphreys:CoxeterGroups}). Denote by $W\mathcal{S}$
the set of all spherical cosets in $W$, that is, the disjoint
union
$$
        W\mathcal{S} = \bigcup_{T \in \mathcal{S}} W/W_T \, .
$$
This is a partially ordered set (\emph{poset}) by inclusion, and
there is a natural $W$-action whose quotient space is
$\mathcal{S}$. The \emph{Davis complex} $\Sigma$ is defined as the
geometric realization $|W\mathcal{S}|$ of the poset
$W\mathcal{S}$, that is, the simplicial complex with one
$n$-simplex for each chain of length $n$
\begin{equation}\label{eqn:ChainDavisComplex}
    w_1W_{T_1} \subset \ldots \subset w_nW_{T_n}\, , \quad T_i \in
    \mathcal{S}\, ,
\end{equation}
and obvious identifications. The $W$-action on $W\mathcal{S}$
induces an action on $\Sigma$ for which this is a proper space:
the stabilizer of the simplex corresponding to the chain
(\ref{eqn:ChainDavisComplex}) is $w_1W_{T_1}w_1^{-1}$, finite.
Moreover, $\Sigma$ admits a CAT(0)-metric \cite{MoussongPhD} and
therefore it is a model of $\underline{E}W$, by the following
result.

\begin{proposition}
If a finite group $H$ acts by isometries on a CAT(0)-space $X$,
then the fixed point subspace $X^H$ is contractible.
\end{proposition}
\begin{proof}
The idea is to find a fixed point for $H$ and then use that if $x,
y \in X$ are fixed by $H$ then so is the geodesic between $x$ and
$y$ (from which the contractibility of $X^H$ follows). The
existence of a fixed point for $H$ is a consequence of the
Bruhat-Tits fixed point theorem \cite[p.~157]{Brown:Buildings}.
\end{proof}

\begin{figure}[htb]
\begin{center}
\includegraphics[scale = 0.8]{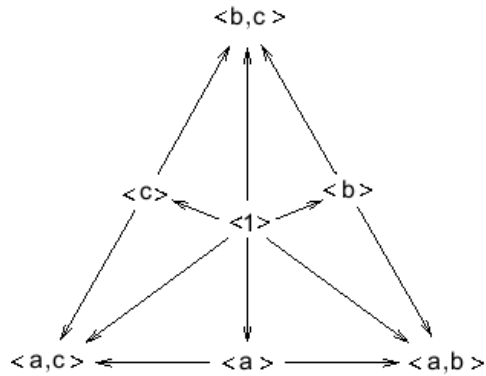}
\caption{Davis complex (quotient) for a triangle group.}
\label{fig:trianglesd}
\end{center}
\end{figure}

\noindent \textbf{Examples:}\begin{itemize} \item[(1)]
\underline{Finite groups.-} \label{example:CoxeterCellFiniteGroup}
Suppose that $W$ is a finite Coxeter group and $C$ is its Coxeter
cell, defined as the convex hull of certain orbit $Wx$. Denote by
$\mathcal{F}(C)$ the poset of faces of the convex polytope $C$.
The correspondence
$$
    w \in W \;\mapsto\; wx \in Wx \  = \textrm{ vertex set of C}
$$
allows us to identify a subset of $W$ with a subset of vertices of
$C$. In fact, a subset of $W$ corresponds to the vertex set of a
face of $C$ if and only if it is a coset of a special subgroup
(\cite[III]{Brown:Buildings}); see Figure \ref{fig:hexagon} for an
example. Hence, we have the isomorphism of posets $W\mathcal{S}
\cong \mathcal{F}(C)$ and the Davis complex $\Sigma =
|W\mathcal{S}|$ is the barycentric subdivision of the Coxeter
cell.

\item[(2)] \underline{Triangle groups.-} These are the Coxeter
groups of rank 3,
$$
    \langle\, a, b, c \;|\; a^2 = b^2 = c^2 = (ab)^p = (bc)^q =
    (c a)^r = 1 \,\rangle\, .
$$
Suppose that $p,q,r \neq \infty$ and $\frac{1}{p} + \frac{1}{q} +
\frac{1}{r} \le 1$. All subsets $T \subsetneq S$ are spherical,
giving the poset (where arrows stand for inclusions) of Figure
\ref{fig:trianglesd}. It can be realized as the barycentric
subdivision of an euclidean or hyperbolic triangle with interior
angles $\frac{\pi}{p}$, $\frac{\pi}{q}$ and $\frac{\pi}{r}$, and
$a$, $b$ and $c$ acting as reflections through the corresponding
sides. The whole model consists on a tessellation of the euclidean
or hyperbolic space by these triangles.
\end{itemize}

\subsection*{A second definition of the Davis complex}\label{section:SecondDefinitionDavisComplex}
There is an alternative description of the Davis complex in terms
of Coxeter cells. Given a poset $\mathcal{P}$ and an element $X
\in \mathcal{P}$, we denote by $\mathcal{P}_{\le X}$ the subposet
of elements in $\mathcal{P}$ less or equal to $X$. Consider a
spherical subset $T \subset S$ and an element $w \in W$.
\begin{proposition}\label{proposition:IsoOfPosets}
There is an isomorphism of posets $\left(W\mathcal{S}\right)_{\le wW_T} \cong
W_T\mathcal{S}$.
\end{proposition}
\begin{proof} Firstly, observe that the poset
$\left(W\mathcal{S}\right)_{\le wW_T}$ is equivalent to
$\left(W\mathcal{S}\right)_{\le W_T}$ via the isomorphism induced
by multiplication by $w^{-1}$. So it suffices to show that
$W\mathcal{S}_{\le W_T} = W_T\mathcal{S}$. A standard element in
the right hand side is a coset $wW_{T'}$ with $w \in W_T$ and $T'
\subset T$ so we have $W_{T'} \subset W_T$ and $wW_{T'} \subset
wW_T = W_T$. On the other hand, if $wW_{T'} \subset W_T$ then $w =
w \cdot 1 \in W_T$. Finally, $T' \subset T$: if $s' \in T'$ then
$ws' \in W_T$ so $s' \in W_T$ and, using the deletion condition,
$s' \in T$.
\end{proof}

\noindent Now, since $W_T$ is a finite Coxeter group, we can
identify the subcomplex $$|\left(W\mathcal{S}\right)_{\le wW_T}| =
|W_T\mathcal{S}|$$ with the barycentric subdivision of the
associated Coxeter cell $C_{W_T}$. From this point of view, the
poset $W\mathcal{S}$ is the union of subposets
$\left(W\mathcal{S}\right)_{\le wW_T}$ with $T$ spherical and $w
\in W/W_T$. That is, $\Sigma$ can be viewed as the union of
Coxeter cells
$$
    \bigcup_{\substack{T \in \mathcal{S} \\ wW_T \in
    W/W_T}} C_{wW_T}
$$
where $C_{wW_T}$ is the copy of $C_{W_T}$ corresponding to the
coset $wW_T$. This union is obviously not disjoint, and the
inclusions and intersections among subposets are precisely the
following.

\begin{lemma} \begin{enumerate}
\item[(i)] $wW_T \subset w'W_{T'}$ if and only if $T \subset T'$
and $w^{-1}w' \in W_{T'}$. \item[(ii)] $wW_T \cap w'W_{T'} =
w_0W_{T \cap T'}$ if there is any $w_0W_T \in wW_T \cap w'W_{T'}$,
empty otherwise.
\end{enumerate}
\end{lemma}

\begin{proof} (i) One implication is obvious, the other uses the
deletion condition in the fashion of the proof of Proposition
\ref{proposition:IsoOfPosets}.

(ii) Straightforward, since $W_T \cap W_{T'} = W_{T \cap T'}\, .$\end{proof}

\noindent Consequently, the intersection of two Coxeter cells is
\begin{equation}\label{eqn:IntersectionCoxCells}
    C_{wW_T} \cap C_{w'W_{T'}} = \left\{ \begin{array}{cl}
        C_{w_0W_{T\cap T'}} & \textrm{if } wW_T \cap w'W_{T'} = w_0 W_{T \cap T'}\\
        \varnothing & \textrm{otherwise} \, .
    \end{array} \right.
\end{equation}

Denote by $\partial C_W$ the \emph{boundary} of the Coxeter cell,
that is, the topological boundary in the ambient space
$\mathbb{R}^{N}$. We have the following explicit description.
\begin{proposition}\label{prop:BoundaryCoxeterCell}
$$\partial C_W = \bigcup_{\substack{T \varsubsetneq
S\\ wW_T \in W/W_T}}wC_{W_T}$$
\end{proposition}
\begin{proof} $C_W$ is a convex polytope whose faces correspond to cosets of
special subgroups (Example on page
\pageref{example:CoxeterCellFiniteGroup}). The barycentric
subdivision $sd(C_W)$ = $|W\mathcal{S}|$ is indeed a cone over its
center (the vertex corresponding to the coset $W_S = W$).
\end{proof}

%%%% NEW SECTION %%%%%%%%%%%%%%%%%%%%%%%%%%%%%%%%%%%%%%%%%%%%%%%%%%%%%
\section{Relative Bredon homology of the Davis complex}\label{section:RelativeCoxeter}
The aim of this section is to deduce a formula for the relative
Bredon homology of the $n$-skeleton of the Davis complex $\Sigma$
with respect to the $(n-1)$-skeleton,
$H^{\mathfrak{Fin}}_*\left(\Sigma_n, \Sigma_{n-1}\right)$, in
terms of the Bredon homology of the Coxeter cells.

We will denote by $\Sigma$ the Davis complex with the Coxeter cell
decomposition explained on page
\pageref{section:SecondDefinitionDavisComplex} and by $\Sigma'$
the original definition as the nerve of $W\mathcal{S}$. Note that
both are $W$-homeomorphic spaces and $\Sigma' = sd(\Sigma)$
(barycentric subdivision) but $\Sigma$ is \emph{not} a
$W$-CW-complex while $\Sigma'$ is (recall that in a $G$-CW-complex
a cell is sent to itself only by the identity map).

\begin{figure}[hbt]
   \begin{center}
   \includegraphics[scale = 0.8]{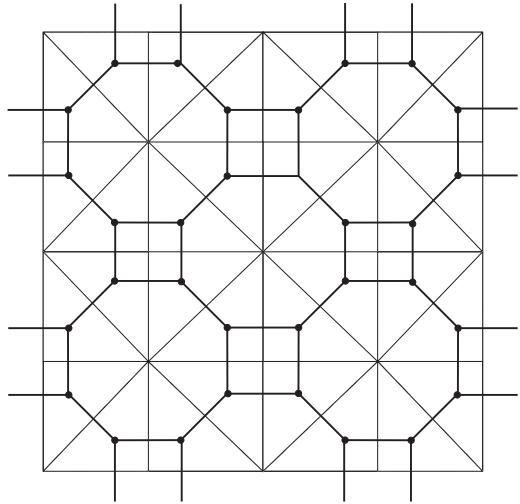}
   \caption[Example of the Coxeter cell structure for $\Delta(2,4,4)$ ]
   {Tessellation of the euclidean plane given by the triangle group $\Delta(2,4,4)$,
   and its dual tessellation.}
   \label{fig:tessellation}
   \end{center}
   \end{figure}

As an example, see Figure \ref{fig:tessellation}: it is the
tessellation of the euclidean plane induced by the triangle group
$\Delta(2,4,4)$ together with the dual tessellation given by
squares and octagons. The latter is $\Sigma$, the Davis complex
given as union of Coxeter cells. The skeleton filtration
correspond to Coxeter cells of rank 0 (points), rank 1 (intervals)
and rank 2 ($2n$-gons).

\begin{defn}
For $n \ge -1$ define $\Sigma_n \subset \Sigma$ as the union
of Coxeter cells corresponding to finite $W_T$ with
$\textrm{rank}(T) \le n$. That is,
$$
    \Sigma_n = \bigcup_{\substack{T \in \mathcal{S},\, |T| \le n\\ wW_T \in
    W/W_T}} C_{wW_T} \, .
$$
\end{defn}
Then $\Sigma_{-1} = \varnothing$, $\Sigma_0$ is a free orbit of
points and $\Sigma_1$ is indeed the Cayley graph of $W$ with
respect to $S$ (see Figure \ref{fig:tessellation} for an example
or compare with Proposition \ref{prop:CoxCellAsInducedSpace}).

We have that $\Sigma_n$ is a $W$-subspace ($w'C_{wW_T} =
C_{w'wW_T}$) and, since the dimension of $C_{W_T}$ is
$\textrm{rank}(T)$, $\Sigma_n$ is indeed the $n$-skeleton of
$\Sigma$.

Our aim is to prove the following theorem.
\begin{theorem}\label{thm:FormulaCoxeterCells}
    For any dimension $n \ge 0$, any degree $i$ and any Bredon
    coefficient system $M$, we have the isomorphism
    $$
        H^{\mathfrak{Fin}}_i \left( \Sigma_n, \Sigma_{n-1}; M
        \right) \cong \bigoplus_{\substack{T \in \mathcal{S}\\ \textup{rank}(T) =
        n}} H^{\mathfrak{Fin}}_i \left( C_{W_T}, \partial
        C_{W_T}; M \right) \, .
    $$
\end{theorem}
\noindent Comments on the statement of the theorem:

(i) We have defined Bredon homology on $G$-CW-complexes so we
implicitly assume that we take the barycentric subdivision on the
spaces.

(ii) The Bredon homology on the left-hand side is with respect to
the $W$-action and the family $\mathfrak{Fin}(W)$. The Bredon
homology on the right-hand side is with respect to the
$W_T$-action and the family $\mathfrak{Fin}(W_T) =
\mathfrak{All}(W_T)$ (all subgroups), for each Coxeter cell.

Firstly, we state some properties of the Coxeter cells that we
will use.
\begin{lemma}\label{lemma:PropertiesCoxCells(3in1)}
\begin{itemize}
    \item[(a)] Suppose $T \in \mathcal{S}$ with $\textrm{rank}(T) \le n$. Then
$$C_{W_T} \cap \Sigma_{n-1} = \left\{ \begin{array}{ll}
            C_{W_T} & \textrm{if } \textrm{rank}(T) < n\, , \\
            \partial C_{W_T} & \textrm{if } \textrm{rank}(T) = n\, .
            \end{array} \right.$$
    \item[(b)] The intersection of any two different Coxeter cells in
$\Sigma_{n}$ is included in $\Sigma_{n-1}$.
    \item[(c)] The Coxeter cell $C_{W_T}$ is a $W_T$--subspace of $\Sigma_n$ and
$\stab_W(x)$ is contained in $W_T$ for all $x \in C_{W_T}$.
\end{itemize}
\end{lemma}

\begin{proof} (a) If $\textrm{rank}(T) < n$ then $C_{W_T}
\subset \Sigma_{n-1}$. If $\textrm{rank}(T) = n$ then (Proposition
\ref{prop:BoundaryCoxeterCell})
$$ \partial C_{W_T} = \bigcup_{\substack{T' \varsubsetneq T\\ wW_{T'} \in W_T/W_{T'}}} wC_{W_{T'}}
\subset \Sigma_{n-1} \, .$$ On the other hand, if $x \in C_{W_T}
\cap \Sigma_{n-1}$ then $x \in C_{wW_{T'}}$ with
$\textrm{rank}(T') \le n-1$ so $x \in C_{W_T} \cap wC_{W_{T'}}$
which, by equation (\ref{eqn:IntersectionCoxCells}), equals to
either empty or $w_0C_{W_{T \cap T'}} \subset \partial C_{W_T}$.

(b) In general, the intersection of two different Coxeter cells of
rank $n$ and $m$ is either empty or another Coxeter cell of rank
strictly less than $\min \{n,m\}$ --- see equation
(\ref{eqn:IntersectionCoxCells}).

(c) The first part is obvious: $wC_{W_T} = C_{wW_T} = C_{W_T}$ if
$w \in W_T$.\\ Suppose now that $x \in C_{W_T}$ and $wx = x$. We
want $w \in W_T$. Identify $sd(C_{W_T}) = |W_T\mathcal{S}|$.
Suppose that $x$ belongs to a cell corresponding to the chain
$W_{T_0} \varsubsetneq \ldots \varsubsetneq W_{T_k}$, $T_k \subset
T$. If $w$ fixes $x$, it has to fix at least one of the vertices
of the cell, i.e., there is an $i$ such that $wW_{T_i} = W_{T_i}$
so that $w \in W_{T_i} \subset W_T$.
\end{proof}

Secondly, we observe that the copies of a Coxeter cell in $\Sigma$
admit an interpretation as induced spaces. If $H$ is a subgroup of
$G$ and $X$ is an $H$-space, the associated \emph{induced}
$G$-space is
$$
    \textup{Ind}_H^G X = G \times_H X\,,
$$
the orbit space for the $H$-action $h \cdot (g,x) = (gh^{-1},hx)$.
The (left) $G$-action on the induced space is given by $g \cdot
\overline{(k,x)} = \overline{(gk,x)}$. This definition carries on
to pairs of $H$-spaces.

\begin{proposition}\label{prop:CoxCellAsInducedSpace} For each spherical $T \subset S$,
there is a $W$-homeomorphism
$$
    \bigcup_{wW_T \in W/W_T} C_{wW_T} \:\cong_W\: W \times_{W_T}
    C_{W_T}\, .
$$
\end{proposition}
\begin{proof}
If $x \in C_{wW_T} = w\cdot C_{W_T}$, write $x = wx_0$ and send
$x$ to $\overline{(w,x_0)} \in W \times_{W_T} C_{W_T}$. This is
well-defined by equation \eqref{eqn:IntersectionCoxCells}:
$x=wx_0=w'y_0$ then there is $w'' \in W_T$ such that $w=w'w''$.
This defines a continuous $W$-map, bijective, with continuous
inverse
\begin{equation*}
    \overline{(w,x)} \mapsto w\cdot x \in C_{wW_T}\, .
    \qedhere
\end{equation*}
\end{proof} \noindent Consequently, we may write $\Sigma_n$, $n
\ge 1$, as a union of induced spaces
$$
    \Sigma_n \cong \bigcup_{|T|\le n} W \times_{W_T}
    C_{W_T}\, .
$$

Next, we will need the following easy consequence of a relative
Mayer-Vietoris sequence.
\begin{proposition}\label{prop:ConsequenceRelativeMayerVietoris}
Let $X$ be a $G$-CW-complex and $Y, A_1, \ldots, A_n$
$G$-sub\-com\-ple\-xes such that $X = A_1 \cup \ldots \cup A_n$
and $A_i \cap A_j \subset Y$ for all $i \neq j$. Write $H_n(-)$
for Bredon homology with some fixed coefficients and with respect
to a family $\F \supset \mathfrak{Iso}(X)$ or, more generally, any
$G$-homology theory. Then
$$
    H_n(X,Y) \cong \bigoplus_{i=1}^n H_n(A_i, A_i \cap Y)\, .
$$
\end{proposition}
\begin{proof}
By induction on $n$. For $n=1$ it is a tautology. For $n>1$ call
$A = A_1$, $B = A_2 \cup \ldots \cup A_n$, $C = A \cap Y$ and $D =
B \cap Y$. The relative Mayer-Vietoris of the CW-pairs $(A,C)$ and
$(B,D)$ is
\begin{equation}\label{eqn:LongExactSeqRelativeM-V}
     \ldots \rightarrow H_i(A \cap B, C \cap D) \rightarrow H_i(A
    , C) \oplus H_i(B, D) \rightarrow H_i(A \cup B, C \cup D) \rightarrow
    \ldots
\end{equation}
Observe that $$\begin{array}{rcl}
  A \cup B &=& X \\
  C \cup D &=& (A \cup B) \cap Y = Y \\
  A \cap B &\subset& Y \\
  C \cap D &=& A \cap B \cap Y = A \cap B
\end{array}$$
Therefore $H_i(A \cap B, C \cap D) = 0$ for all $i$ and the
sequence \eqref{eqn:LongExactSeqRelativeM-V} gives isomorphisms
$$
    H_i(X, Y) \cong H_i(A,C) \oplus H_i(B,D)
$$
Now apply induction to $X'=B$, $Y'=B\cap Y$ and the result
follows.
\end{proof}

Finally, recall that if $h_*^?$ is an \emph{equivariant homology
theory} (see, for instance, \cite[\S{1}]{Luck02}), it has an
\emph{induction structure}. Hence, for a pair of $H$-spaces
$(X,A)$,
$$
    h^G_n(\textup{Ind}_H^G(X,A)) \cong h^H_n(X,A)\, .
$$
Bredon homology has an induction structure \cite{MS95}; in
particular,
$$
    H^{\Fin(W)}_n\left( \textup{Ind}_{W_T}^W ( C_{W_T},
    \partial C_{W_T}) \right) \cong
    H^{\Fin(W_T)}_n \left( C_{W_T}, \partial C_{W_T} \right)\, .
$$
\begin{proof}(\emph{Theorem \ref{thm:FormulaCoxeterCells}})
Write $H_n(-)$ for Bredon homology of proper $G$-spaces with
coefficients in the representation ring. Let $T_1, \ldots , T_m$
be the spherical subsets of generators up to rank $n$. Define
$$
    A_i = W \times_{W_{T_i}} C_{W_{T_i}} \quad 0 \le i \le m\, .
$$
Then $\Sigma_n = \bigcup A_i$ and $A_i \cap A_j \subset
\Sigma_{n-1}$ for all $i\neq j$, by Lemma
\ref{lemma:PropertiesCoxCells(3in1)}(b). By Proposition
\ref{prop:ConsequenceRelativeMayerVietoris}
$$
    H_n(\Sigma_n,\Sigma_{n-1}) \cong \bigoplus_{i=1}^n H_n(A_i, A_i \cap \Sigma_{n-1})\, .
$$
If $\textup{rank}(T_i) < n$ then $A_i \subset \Sigma_{n-1}$ and
the corresponding term is zero. If $\textup{rank}(T_i) = n$ then
$A_i \cap \Sigma_{n-1} = \partial A_i$ (Lemma
\ref{lemma:PropertiesCoxCells(3in1)}(a)) and, by the induction
structure,
\begin{equation*}
    H_n(A_i, \partial A_i) \cong H_n(C_{W_{T_i}}, \partial
    C_{W_{T_i}})\,. \qedhere
\end{equation*}
\end{proof}\label{addCorollaries}
\begin{corollary}
Write $H_i(-)$ for Bredon homology with respect to finite
subgroups and some fixed Bredon coefficient system. For each $n
\ge 1$, there is a long exact sequence
$$
     \ldots \rightarrow H_i^{\mathfrak{Fin}}(\Sigma_{n-1}) \rightarrow
     H_i^{\mathfrak{Fin}}(\Sigma_n) \rightarrow
     \bigoplus_{\substack{T \in \mathcal{S}\\ \textup{rank}(T) =n}}
     H^{\mathfrak{Fin}}_i \left( C_{W_T}, \partial
        C_{W_T} \right) \rightarrow
    \ldots
$$
\end{corollary}

\begin{remark}
Theorem \ref{thm:FormulaCoxeterCells} also holds for any (proper)
equivariant homology theory, in particular, for equivariant
$K$-homology.
\end{remark}

\section{Relative Bredon homology of some Coxeter cells}
\begin{figure}[hbt]
   \begin{center}
   \includegraphics[scale = 0.4]{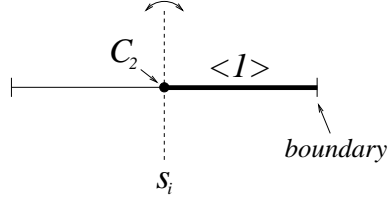}
   \caption{Coxeter cell of rank 1}
   \label{fig:interval}
   \end{center}
   \end{figure}

\begin{figure}[hbt]
   \begin{center}
   \includegraphics[scale = 0.35]{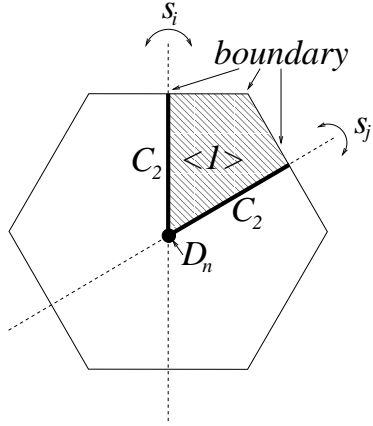}
   \caption{Coxeter cell of rank 2}
   \label{fig:hexagon4}
   \end{center}
   \end{figure}

As an application of Theorem \ref{thm:FormulaCoxeterCells}, we
compute the relative Bredon homology of $(\Sigma_n,
\Sigma_{n-1})$, with coefficients in the representation ring, for
the first cases $n = 0, 1, 2$. To do so, we need the Bredon
homology of $(C_{W_T}, \partial C_{W_T})$ when
$T \subset S$ spherical with $\textrm{rank}(T) = 0,1,2$.\\[10pt]
\framebox{$\textrm{rank}(T)$ = 0} $T = \varnothing$, $W_T = \{1\}$
and $C_{W_T}$ is a point. So the homology is $R_\mathbb{C}(\{1\})
\cong \mathbb{Z}$ at degree 0 and vanishes
elsewhere.\\[10pt]
\framebox{$\textrm{rank}(T)$ = 1} $T = \{s_i\}$, $W_T$ is cyclic of order two and $C_{W_T}$ is
an interval. For the relative homology, we consider one 0-cell
and one 1-cell, with stabilizers cyclic order two and trivial
respectively (see Figure \ref{fig:interval}).%
   The associated Bredon chain complex
$$ \begin{array}{ccccccc}
    0 &\longrightarrow &R_\mathbb{C}(\{1\}) &
    \stackrel{\Ind}{\longrightarrow} & R_\mathbb{C}\left(
    C_2 \right) &\longrightarrow &0\\
     & & 1 & \mapsto & (1,1) & &
    \end{array}
$$
gives
\begin{equation}\label{eqn:RelBredonHomCoxCellRank2}
    H_i \left( C_{W_T}, \partial C_{W_T}; \RR \right) = \left\{
            \begin{array}{cc}
                \mathbb{Z} & i = 0\\
                0 & i \not = 0
            \end{array} \right. \, .
\end{equation}\\[6pt]
\framebox{$\textrm{rank}(T)$ = 2} $T = \{s_i, s_j\}$, $W_T$ is
dihedral of order $2n$ ($n = m_{ij} \not = \infty$)  and $C_{W_T}$
is a $2n$-gon. The orbit space is a sector of the polygon (see
Figure \ref{fig:hexagon4}) and the relative Bredon chain complex
is
$$
    0 \longrightarrow R_\mathbb{C}\left( \{1\} \right)
    \stackrel{d_2}{\longrightarrow}
    R_\mathbb{C}\left( C_2 \right) \oplus R_\mathbb{C} \left( C_2 \right)
    \stackrel{d_1}{\longrightarrow}
    R_\mathbb{C}\left( D_n \right)
    \longrightarrow 0\, ,
$$
where the differentials are given by induction among
representation rings. Inducing the trivial representation always
yields the regular one, so $d_2$ in matrix form is
$$
        d_2 \equiv \left( \begin{array}{cccc}
                            1 & 1 & 1 & 1 \end{array} \right)
        \sim        \left( \begin{array}{cccc}
                            1 & 0 & 0 & 0 \end{array} \right)
                    \, ,
$$
the right-hand side being its Smith normal form. For $d_1$, recall
the character tables of the cyclic group of order two and dihedral
group (with Coxeter generators $a$ and $b$):

\begin{flushleft}
\begin{tabular}[t]{ll}
\begin{tabular}[t]{c|rr}
       $C_2$ & $1$ & $a$ \\
       \hline
       $\rho_1$ & $1$ & $1$ \\
       $\rho_2$ & $1$ & $-1$ \\
\end{tabular}
&
\begin{tabular}[t]{c|cc}
       $D_n$ & $(ab)^k$ & $b(ab)^k$ \\
       \hline
       $\chi_1$ & $1$ & $1$ \\
       $\chi_2$ & $1$ & $-1$ \\
       $\widehat{\chi_3}$ & $(-1)^k$ & $(-1)^k$ \\
       $\widehat{\chi_4}$ & $(-1)^k$ & $(-1)^{k+1}$ \\
       $\phi_l$ & $2\cos\left( 2\pi l k/n\right)$ & $0$ \\
\end{tabular}
\end{tabular}
\end{flushleft}
where $0 \le k \le n-1$, $l$ varies from 1 to $n/2$ ($n$ even) or
$(n-1)/2$ ($n$ odd) and the hat \ $\widehat{ }$ \ denotes a
character which only appears when $n$ is even. The induced
representations corresponding to the inclusions $\langle a \rangle
\subset \langle a,b \rangle$ and $\langle b \rangle \subset
\langle a,b \rangle$ are
\begin{eqnarray*}
    \textrm{Ind}\,(\rho_1) &=& \chi_1 + \widehat{\chi_3} + \sum \phi_l\, ,\\
    \textrm{Ind}\,(\rho_2) &=& \chi_2 + \widehat{\chi_4} + \sum
    \phi_l\ .
\end{eqnarray*}
Consequently, the differential $d_1$ is written in matrix form as
$$
     d_1 \equiv \left(
            \begin{array}{ccccccc}
                1 & 0 & \widehat{1} & \widehat{0} & 1 & \ldots & 1 \\
                0 & 1 & \widehat{0} & \widehat{1} & 1 & \ldots & 1 \\
                1 & 0 & \widehat{0} & \widehat{1} & 1 & \ldots & 1 \\
                0 & 1 & \widehat{1} & \widehat{0} & 1 & \ldots & 1 \\
            \end{array} \right) \, ,
$$
which reduces by elementary operations to
$$
            \left( \begin{array}{c|c} I_3 & 0 \\ \hline 0 & 0 \end{array}\right) \ \textrm{ if $n$ is
            even}\, , \quad
            \left( \begin{array}{c|c} I_2 & 0\\ \hline 0 & 0 \end{array}\right) \ \textrm{ if $n$ is odd}\, . \\
$$
Therefore, the relative Bredon homology of a Coxeter cell of rank
two and $n = m_{ij} \not = \infty$ is
$$ \begin{array}{llll}
        \textrm{$n$ even:} & H_0 \,=\, \mathbb{Z}^{c(D_n)-3} =
        \mathbb{Z}^{n/2}\, , & &  H_i = 0 \ \;\forall \ i \ge 1 \, ;\\
        \textrm{$n$ odd:} & H_0 \,=\, \mathbb{Z}^{c(D_n)-2} =
        \mathbb{Z}^{(n-1)/2}\, , & H_1 = \mathbb{Z}\,, & H_i = 0
        \ \;\forall \ i \ge 2 \, .
    \end{array}
$$
Here we write $c(H)$ for the number of conjugacy classes in a
finite group $H$, and $\Z^n$, or sometimes $n \cdot \Z$, for
$\oplus_{i=1}^n \Z$.

We can now use Theorem \ref{thm:FormulaCoxeterCells} to deduce the
lower-rank relative Bredon homology of $\Sigma$. We omit the
family $\mathfrak{Fin}$ and coefficients $\RR$ for clarity.

\begin{proposition}\label{prop:LowerRankRelativeHomDavisCx}
Let $(W,S)$ be a Coxeter system and $\Sigma$ its associated
Coxeter complex. Then
$$
\begin{array}{rclc}
        H_i \left( \Sigma_0, \Sigma_{-1} \right) & = & \left\{
            \begin{array}{cl} \mathbb{Z} & i = 0\\ 0 & i \neq
            0 \end{array} \right. \, .\\[1.3em]
        H_i \left( \Sigma_1, \Sigma_0 \right) & = & \left\{
            \begin{array}{cl} |S| \cdot \mathbb{Z} & i = 0\\ 0 & i \neq
            0 \end{array} \right. \, .
    \end{array}
$$
$$
        H_0 \left( \Sigma_2, \Sigma_1 \right)  =
            \frac{1}{2}\left( \sum_{\textrm{even, }i <
            j} m_{ij} + \sum_{\textrm{odd, }\\
            i<j} (m_{ij}-1) \right) \cdot
            \mathbb{Z} \, ,
$$
$$
        H_1 \left( \Sigma_2, \Sigma_1 \right)  = |\{ m_{ij} \textrm{ odd}, i<j \}| \cdot
        \mathbb{Z} \, ,\quad
       H_i \left( \Sigma_2, \Sigma_1 \right)  = 0  \ \;(i \ge
       2)\, .
$$
\end{proposition}

%%%% NEW SECTION %%%%%%%%%%%%%%%%%%%%%%%%%%%%%%%%%%%%%%%%%%%%%%%%%%%%%
\section[Bredon homology of right-angled Coxeter groups]
{Bredon homology of the right-angled Coxeter
groups}\label{section:RightAngled} A Coxeter system $(W,S)$ is
\emph{right-angled} if $m_{ij} \in \{ 2, \infty \}$ for all $i,
j$. Thus, the only spherical subgroups are direct products of
cyclic groups of order 2.

Recall that $C_{W_1 \times W_2} = C_{W_1} \times C_{W_2}$. Hence,
if $W_T$ is a product of cyclic groups of order 2, we can use a
relative K\"unneth formula (Section \ref{section:KunnethFormulas})
and equation (\ref{eqn:RelBredonHomCoxCellRank2}) to deduce
\begin{equation}\label{eqn:HomCoxCellRA}
    H_i \left( C_{W_T}, \partial C_{W_T} \right) = \left\{
            \begin{array}{cc}
                \mathbb{Z} & i = 0\, ,\\
                0 & i \not = 0\, .
            \end{array} \right.
\end{equation}
Call $s(n)$ the number of distinct spherical subsets $T \subset S$
of rank $n$ ($n \ge 0$). For instance, $s(0) = 1$, $s(1) = |S|$
and $s(2) = |\{m_{ij} \neq \infty,\, i\!<\!j \}|$. By Theorem
\ref{thm:FormulaCoxeterCells} and (\ref{eqn:HomCoxCellRA})
$$
    H_i \left( \Sigma_n, \Sigma_{n-1} \right) =
    \left\{ \begin{array}{cl} s(n) \cdot \mathbb{Z} & i = 0\, ,\\
    0 & i \neq 0 \, ,\end{array} \right.
$$
for all $n \ge 0$. The long exact sequence of the pair $(\Sigma_n,
\Sigma_{n-1})$ gives $H_i(\Sigma_n) \cong H_i(\Sigma_{n-1})$ for
$i > 0$ and the split exact sequence
$$
    0 \longrightarrow H_0(\Sigma_{n-1}) \longrightarrow H_0(\Sigma_{n})\longrightarrow H_0(\Sigma_n,
    \Sigma_{n-1}) \longrightarrow 0 \, .
$$
Using induction on $n$,
$$  \begin{array}{l}
    H_i(\Sigma_n) = H_i(\Sigma_0) = 0 \quad \forall \ i>0\,
    ,\\[1ex]
    H_0(\Sigma_n) = H_0(\Sigma_{n-1}) \oplus H_0(\Sigma_n,
    \Sigma_{n-1}) = \left( s(0) + \ldots + s(n) \right) \cdot
    \mathbb{Z} \, .
    \end{array}
$$
This gives the Bredon homology of $(W,S)$ when $n = |S|$.
\begin{theorem} \label{thm:BredonHomRA}
The Bredon homology of a right-angled Coxeter
group $(W,S)$ with respect to the family $\mathfrak{Fin}(W)$
and coefficients in the representation ring is
$$
    H^\mathfrak{Fin}_i \left( W; \RR
    \right) = \left\{ \begin{array}{ll} s
    \cdot \mathbb{Z} & i = 0\\
    0 & i \neq 0
    \end{array} \right.
$$
where $s = s(0) + \ldots + s(|S|)$, the number of spherical
subsets $T \subset S$.
\end{theorem}

Finally, since the Bredon homology concentrates at degree 0, it
coincides with the equivariant $K$-homology (Proposition
\ref{prop:LowDegreeBredonHomology}). Also, the Baum-Connes
conjecture holds for Coxeter groups, so we have the following
corollary.

\begin{corollary}\label{cor:KhomologyRA}
If $W$ is a right-angled Coxeter group, the equivariant
$K$-homo\-lo\-gy of $\underline{E}W$ coincides with its Bredon
homology at degree 0 and 1 respectively (given in the previous
theorem). This corresponds to the topological $K$-theory of the
reduced $C^*$-algebra of $W$, via the Baum-Connes assembly map.
\end{corollary}

%%%% NEW SECTION %%%%%%%%%%%%%%%%%%%%%%%%%%%%%%%%%%%%%%%%%%%%%%%%%%%%%
\section{Bredon homology of even Coxeter groups}\label{section:EvenCoxeter}
A Coxeter system $(W,S)$ is \emph{even} if $m_{ij}$ is even or
infinite for all $i,j$. Right-angled Coxeter groups are then a
particular case. Again, the only spherical subgroups are direct
products of cyclic (order two) and dihedral, since any irreducible
finite Coxeter group with more than two generators have at least
one odd $m_{ij}$.

Suppose that $W_T$ = $D_{m_1} \times \ldots \times D_{m_r} \times
(C_2)^k$, with $m_i$ ($1 \le i \le r$) even numbers and $k \ge 0$.
We know that the relative Bredon homology of the Coxeter cells
corresponding to dihedral $D_{m_i}$ and $C_2$ concentrates at
degree 0 and
$$
    H_0 = (m_i/2) \cdot \mathbb{Z} \, ,\quad \textrm{respectively} \quad H_0 =
    \mathbb{Z} \, .
$$
By the relative K\"unneth formula we have
$$
    H_i \left( C_{W_T}, \partial C_{W_T} \right) = \left\{
    \begin{array}{cl}
        \left( m_1 \cdot \ldots \cdot m_r \big/ 2^r\right) \cdot \mathbb{Z} &
        i=0\\ 0 & i\neq 0 \end{array} \right. \, .
$$
Alternatively, if $T = \{ t_1, \ldots , t_n \}$ is a spherical
subset of rank $n$,
\begin{equation}\label{eqn:RelativeBredonHomEvenCoxGrps}
    H_0 \left( C_{W_T}, \partial C_{W_T} \right) = \prod_{i<j} m_{ij}/2 \cdot \Z\, ,
\end{equation}
and 0 at any other degree. This expression is valid for any $n
\ge 0$, with the convention that an empty product equals 1.

Now we can compute the relative homology groups of the skeleton
filtration of $\Sigma$. They concentrate at degree zero, since so
does the homology of the Coxeter cells. We already know that
$$
    H_0 \left( \Sigma_0, \Sigma_{-1} \right) = \mathbb{Z},
    \quad
    H_0 \left( \Sigma_1, \Sigma_0 \right) = |S| \cdot
    \mathbb{Z},
    \quad
    H_0 \left( \Sigma_2, \Sigma_1 \right) = \sum_{\substack{i<j\\ m_{ij} \neq \infty}} m_{ij}/2 \ \cdot
    \mathbb{Z}\; .
$$
In general, using Theorem \ref{thm:FormulaCoxeterCells} and
equation (\ref{eqn:RelativeBredonHomEvenCoxGrps}), we have that
$H_0 \left( \Sigma_n, \Sigma_{n-1} \right)$, $n \ge 0$ is a free
abelian group of rank
$$
    r(n) = \sum_{\substack{T \in \mathcal{S}\\ \textup{rank}(T)=n}} \prod_{\substack{s_i, s_j \in
    T\\i<j}} m_{ij}/2 \, ,
$$
with the convention that an empty sum is 0 and an empty
product is 1. As before, the long exact sequences of
$(\Sigma_n, \Sigma_{n-1})$ and induction on $n$ gives
$$ \begin{array}{l}
    H_i \left(\Sigma_n \right) = H_i \left(\Sigma_0 \right) =
    0 \quad \forall \; i \ge 1 \ \textrm{ and}\\
    H_0 \left(\Sigma_n \right) = H_0 \left(\Sigma_{n-1}
    \right) \oplus H_0 \left(\Sigma_n, \Sigma_{n-1} \right) \,
    .
\end{array}
$$
Define
$$
    r = r(0) + \ldots + r(|S|) = \sum_{T \in \mathcal{S}}
    \prod_{\substack{s_i, s_j \in
    T\\i<j}} m_{ij}/2 \, .
$$
\begin{theorem}\label{thm:BredonHomEven}
    The Bredon homology of an even Coxeter group is the free
    abelian group of rank $r$ at degree 0
    and vanishes at any other degree.
\end{theorem}

\begin{remark} This theorem yields Theorem \ref{thm:BredonHomRA}
for a right-angled Coxeter group. It agrees as well with the
Bredon homology of the infinite dihedral group and the triangle
groups with $p$, $q$ and $r$ even (Section
\ref{section:LowRankCoxeterGroups}).
\end{remark}

As before, the Bredon homology coincide with the equivariant
$K$-homology.

\begin{corollary}\label{cor:KhomologyEven}
If $W$ is an even Coxeter group, the equivariant $K$-homology of
$\underline{E}W$ coincides with its Bredon homology at degree 0
and 1 respectively (given in the previous theorem). This
corresponds to the topological $K$-theory of the reduced
$C^*$-algebra of $W$ via the Baum-Connes assembly map.
\end{corollary}

%%%% NEW SECTION %%%%%%%%%%%%%%%%%%%%%%%%%%%%%%%%%%%%%%%%%%%%%%%%%%%%%
\section{Bredon homology of low-rank Coxeter groups}\label{section:LowRankCoxeterGroups}
We briefly recall, for completeness, the Bredon homology of the
Coxeter groups of rank up to three. If a group $H$ is finite, a
one-point space is a model of \EG{} so its Bredon homology reduces
to the representation ring $R_\mathbb{C}(H)$ at degree 0 and
vanishes elsewhere, and so its equivariant $K$-homology. The first
infinite Coxeter groups are the infinite dihedral
$$
    D_\infty = \left< a, b \,|\, a^2 = b^2 = 1 \right> \, ,
$$
and the euclidean and hyperbolic triangle groups
$$
    \left< a, b, c \,|\, a^2 = b ^2 = c^2 = 1,\
        (ab)^p = (bc)^q = (ca)^r = 1 \right> \, ,
$$
where $2 \leq p, q, r \le \infty$ and $1/p + 1/q + 1/r \leq 1$.
All cases can be done by direct computation from the correponding
orbit spaces $\underline{B}G$ (\cite[Section
5.2]{SanchezGarciaPhD}). The resulting Bredon homology
concentrates at degree $i \le 1$ so it coincides with their
$K^G$-homology. The results are the following (write $H_i(G)$ for
$H^\Fin_i\left(\EG; \RR\right)$, $n \cdot \Z$ for $\oplus_{i=1}^n
\Z$, and $c(H)$ for the number of conjugacy classes on a finite
group $H$).
$$\begin{array}{rcl}\label{eqn:BredonHomLowerRankCoxGrps}
H_i\left(D_\infty\right) &=& \left\{ \begin{array}{ll}
                            3 \cdot \Z & i = 0 \\ 0 & i \not = 0
                            \end{array} \right. \\[2em]
H_i\left(\Delta(\infty, \infty, \infty)\right) &=& \left\{
                            \begin{array}{ll}
                            4 \cdot \Z & i = 0 \\ 0 & i \not = 0
                            \end{array} \right. \\[2em]
H_i\left(\Delta(p, \infty, \infty)\right) &=& \left\{
                            \begin{array}{ll}
                            \left(c(D_p)+1\right) \cdot \Z & i = 0 \\ 0 & i \not = 0
                            \end{array} \right. \\[2em]
H_i\left(\Delta(p, q, \infty)\right) &=& \left\{
                            \begin{array}{ll}
                            \left(c(D_p)+c(D_q)-2\right) \cdot \Z & i = 0 \\ 0 & i \not = 0
                            \end{array} \right. \\[2em]
H_0\left(\Delta(p, q, r)    \right) &=& \left\{\begin{array}{ll}
                            \left(c(D_p)+c(D_q)+c(D_r)-4\right) \cdot \Z & p, q, r
                            \textrm{ odd}\\
                            \left(c(D_p)+c(D_q)+c(D_r)-5\right) \cdot \Z & \textrm{ otherwise}\\
                            \end{array} \right. \\[2em]
H_1\left(\Delta(p, q, r)\right) &=& \left\{
                            \begin{array}{ll}
                            \Z & p, q \textrm{ and } r \textrm{ odd}\\
                            0 & \textrm{ otherwise}
                            \end{array} \right. \\[2em]
H_i\left(\Delta(p, q, r)\right) &=& 0 \qquad \textrm{for } i \ge
2\, .
\end{array}
$$
(Recall that $c(D_n)$ is $n/2 +3$ if $n$ if even and $(n-1)/2 + 2$
if $n$ is odd.)
 \begin{remark}
    L\"uck and Stamm \cite{LS00} studied the equivariant
    $K$-homology of cocompact (= compact quotient) planar groups.
    Our groups above, except those triangle groups with $p$, $q$ or $r$
    being $\infty$, are examples of cocompact planar groups and
    our results also follow from Part (c) of Theorem 4.31 in their
    article.
 \end{remark}
 Finally, observe that the results above agree with Theorems
 \ref{thm:BredonHomRA} and \ref{thm:BredonHomEven} when $p$, $q$
 and $r$ are even or infinity.

%%%% NEW SECTION %%%%%%%%%%%%%%%%%%%%%%%%%%%%%%%%%%%%%%%%%%%%%%%%%%%%%
\section{A GAP routine for Coxeter groups}
We have implemented a GAP \cite{GAP} program to obtain the Bredon
homology with coefficients in the representation ring of virtually
any Coxeter group, from its Coxeter matrix.

The GAP routine works as follows. Firstly, we generate the
cellular decomposition, cell stabilizers and boundary of the
corresponding Davis complex. To do that, we have written functions
to find the irreducible components of a Coxeter group and to
decide whether the group it is finite (comparing each irreducible
component with the finite Coxeter groups of the same rank). We
generate a multidimensional list \verb+CHAINS+ such that
\verb+CHAINS[i]+ is the list of all chains $T_1 < \ldots < T_i$ of
length $i$ ($<$ means strict inclusion) of subsets of $\{1,
\ldots, N\}$ where $N = \textrm{rank}(W,S)$. Note that a cell
$\sigma$ in the quotient space $\Sigma/W$ corresponding to $\{T_1
< \ldots < T_n \}$ has dimension $n-1$, stabilizer $W_{T_1}$ and
boundary
$$
        \partial \sigma = \sum_{k=1}^{n}(-1)^k \{T_1 < \ldots <
        \widehat{T_i} < \ldots T_n\}\, .
$$
Secondly, we use this information as initial data to another
procedure which computes the Bredon homology with coefficients in
the representation ring of a finite proper $G$-CW-complex. Note
that, however, the computing time increases exponentially with the
number of generators.

%%% THE BIBLIOGRAPHY %%%%%%%%%%%%%%%%%%%%%%%%%%%%%%%%%%%%%%%%%%%%
\bibliographystyle{amsplain}
\bibliography{MyBibliography}

\providecommand{\bysame}{\leavevmode\hbox to3em{\hrulefill}\thinspace}
\providecommand{\MR}{\relax\ifhmode\unskip\space\fi MR }
% \MRhref is called by the amsart/book/proc definition of \MR.
\providecommand{\MRhref}[2]{%
  \href{http://www.ams.org/mathscinet-getitem?mr=#1}{#2}
}
\providecommand{\href}[2]{#2}
\begin{thebibliography}{10}

\bibitem{BCH94}
P.~Baum, A.~Connes, and N.~Higson, \emph{Classifying spaces for proper actions
  and {$K$-theory} of group {$C^*$-algebras}}, Contemporary Mathematics
  \textbf{167} (1994), 241--291.

\bibitem{BJS88}
M.~Bo{\.z}ejko, T.~Januszkiewicz, and R.~J. Spatzier, \emph{Infinite {C}oxeter
  groups do not have {K}azhdan's property}, J. Operator Theory \textbf{19}
  (1988), no.~1, 63--67.

\bibitem{Brown:Buildings}
K.~S. Brown, \emph{Buildings}, Springer-Verlag, Berlin, 1988.

\bibitem{DL98}
J.~F. Davis and W.~L{\"u}ck, \emph{Spaces over a {Category} and {Assembly}
  {Maps} in {Isomorphism} {Conjectures} in {$K$}- and {$L$}-{Theory}},
  $K$-Theory \textbf{15} (1998), 201--252.

\bibitem{Davis83}
M.~W. Davis, \emph{Groups generated by reflections}, Ann. Math. \textbf{117}
  (1983), 293--324.

\bibitem{GAP}
The GAP~Group, \emph{{GAP -- Groups, Algorithms, and Programming, Version
  4.4}}, 2004, \verb+(http://www.gap-system.org)+.

\bibitem{GG99}
M.~Golasi{\'n}ski and D.~L. Gon{\c{c}}alves, \emph{Generalized
  {E}ilenberg-{Z}ilber type theorem and its equivariant applications}, Bull.
  Sci. Math. \textbf{123} (1999), no.~4, 285--298.

\bibitem{HK97}
N.~Higson and G.~Kasparov, \emph{Operator {$K$}-theory for groups which act
  properly and isometrically on {H}ilbert space}, Electron. Res. Announc. Amer.
  Math. Soc. \textbf{3} (1997), 131--142 (electronic).

\bibitem{Humphreys:CoxeterGroups}
J.~E. Humphreys, \emph{Reflection groups and {C}oxeter groups}, Cambridge
  Studies in Advanced Mathematics, vol.~29, Cambridge University Press, 1990.

\bibitem{Joa03}
M.~Joachim, \emph{{$K$}-homology of {$C\sp \ast$}-categories and symmetric
  spectra representing {$K$}-homology}, Math. Ann. \textbf{327} (2003), no.~4,
  641--670.

\bibitem{LeonardiPhD_inprep}
F.~Leonardi, \emph{\textup{Doctoral Thesis}}, ETH Zurich (in preparation).

\bibitem{Luck02}
W.~L{\"u}ck, \emph{Chern characters for proper equivariant homology theories
  and applications to {$K$}- and {$L$}-theory}, J. Reine Angew. Math.
  \textbf{543} (2002), 193--234.

\bibitem{LS00}
W.~L{\"u}ck and R.~Stamm, \emph{Computations of {$K$} and {$L$}-theory of
  cocompact planar groups}, $K$-Theory \textbf{21} (2000), 249--292.

\bibitem{MislinValette:ProperGroupActions}
G.~Mislin and A.~Valette, \emph{Proper group actions and the {B}aum-{C}onnes
  conjecture}, Advanced Courses in Mathematics. CRM Barcelona, Birkh\"auser
  Verlag, 2003.

\bibitem{MS95}
I.~Moerdijk and J.-A. Svensson, \emph{A {S}hapiro lemma for diagrams of spaces
  with applications to equivariant topology}, Compositio Math. \textbf{96}
  (1995), no.~3, 249--282.

\bibitem{MoussongPhD}
G.~Moussong, \emph{Hyperbolic {Coxeter} groups}, Doctoral {Thesis}, Ohio State
  University, 1988.

\bibitem{SanchezGarciaPhD}
R.~S\'anchez-Garc\'\i{}a, \emph{Equivariant {$K$-homology} of the classifiying
  space for proper actions}, Ph.D. thesis, University of Southampton, 2005.

\end{thebibliography}

%%%% END OF DOCUMENT %%%%%%%%%%%%%%%%%%%%%%%%%%%%%%%%%%%%%%%%%%%%
\end{document}